\documentstyle[12pt]{article}
\begin{document}
\title{\bf\huge A one-dimensional many-body 
integrable model from $Z_n$ Belavin model with open boundary conditions
}

\author{
{\bf
Heng Fan$^{a,b}$, Bo-Yu Hou$^b$, Guang-Liang Li $^b$}, 
\\
{\bf Kang-Jie Shi$^{a,b}$,  Yan-Shen Wang$^{a,c}$} \\
$^{a}$ CCAST(World Laboratory),\\ 
P. O. Box 8730, Beijing 100080, P. R. China \\
$^{b}$ Institute of Morden Physics, Northwest University \\
P. O. Box 105, Xi'an, 710069, P. R. China\\
$^{c}$ Zhejiang Institute of Morden Physics, Zhejiang University \\
Hangzhou, 310027, P. R. China}
\maketitle

\begin{abstract}
We use factorized $L$ operator to construct an 
integrable model with open boundary conditions.
By taking trigonometic limit($\tau \rightarrow \sqrt{-1}\infty$) and
scaling limit($\omega \rightarrow 0$), we get a Hamiltonian of
a classical integrable system. 
It shows that this integrable system is similar
to those found by Calogero et al.
\end{abstract}
\newpage
\renewcommand{\thesection}{\Roman{section}}
\baselineskip 1truecm

\newpage
\section{Introdction}
In the last decades, 
a series of one-dimensional integrable
many-body systems have been found [1-6]. 
One effective method to prove the integrability of the
many-body systems is the Lax representation,
which means that we can construct the complete set
of integrals of the motion.
The Lax representation for 
the elliptic Calogero-Moser model was found by 
Krichever [7]. And the Ruijsenaars-Macdonald's 
commuting difference operators can also prove
the integrability of the many-body systems[5,6,8,9].

It is known that there are a lot of two-dimensional
exactly sovled models in
statistical mechanics[10-15]. The integrability is
proved by the commuting transfer matrix.
Recently, some relations between the one-dimensional
many-body systems and the two-dimensional sovable models
are found, see references [16-19] and the references therein.
Hasegawa [16] found that the $L$ operator [20-22] for
the two-dimensional $Z_n$ Belavin
model can be related to the Krichever's Lax matrix [7].
And the commuting difference
operators given by the $L$ operators 
is similar to the Ruijsenaars-Macdonald's
difference operators. Hasegawa studied the $Z_n$ Belavin model
by imposing the periodic boundary conditions. We will study
the $Z_n$ Belavin model by imposing the open boundary conditions,
for the case of open boundary conditions see references [23,24]
and the references therein.
Using the factorized $L$ operators and one solution of the
reflection equation, we can constrcut a cummuting difference
operators which should be equivalent to the Ruijsenaars-Macdonald's
difference operators. By taking a special limit, we find a 
trigonometic integrable model which is similar to the trigonometric
model found by Calogero et al [2-4].
Principally,
we can obtain a classical integrable system if we can find a 
solution to the reflection equation.

The paper is organized as follows: In section 2, we will introduce the
$Z_n$ symmetric Belavin model. The factorized $L$ operator will 
be given in section 3. The commuting difference operators connected
with the transfer matrix with open boundary conditions will be
given in section 4. The special limit is taken in section 5,
the integrable model is found in this section. Section 6 containes
a summary and discussions.

\section{The $Z_n$ symmetric Belavin R matrix}
The $Z_n$ symmetric Belavin R matrix [14,15] is given as
\begin{equation}
R_{12}= \frac{1}{n} \sum_{\alpha \in Z_n^2}
W_{\alpha}(z)I_{\alpha } \otimes I_{\alpha}^{-1},
\end{equation}
with
$$ W_{\alpha}(z) = \frac{\sigma_{\alpha}(z+\eta)}{\sigma_{\alpha}(\eta)}, 
I_{\alpha}=g^{\alpha_2}h^{\alpha_1}, h_{ij}=\delta_{i+1,j}, $$

$$ g_{ij} = \omega^i \delta_{i,j}, 
\omega=e^{\frac{2\pi\sqrt{-1}}{n}}, (i,j \in Z_n),$$

$$ \sigma_{\alpha}(z) \equiv \theta \left [\begin{array}{c}
      \frac{1}{2}+\frac{\alpha_1}{n}\\
      \frac{1}{2}+\frac{\alpha_2}{n}
  \end{array}\right ](z,\tau),$$

 \begin{equation} \theta \left [\begin{array}{c} a\\ b
  \end{array}\right ](z,\tau) 
  \equiv \sum_{m \in Z}e^{\pi\sqrt{-1}
  (m+a)^2\tau+2\pi \sqrt{-1} (m+a)(z+b)}.
  \end{equation}

This R-matrix satisfy the Yang-Baxter equation(YBE) [10,11],
   \begin{equation}
   R_{12}(z_1-z_2)R_{13}(z_1-z_3)R_{23}(z_2-z_3) 
   =R_{23}(z_2-z_3)R_{13}(z_1-z_3)R_{12}(z_1-z_2),
   \end{equation}
where $ R_{12}(z), R_{13}(z) $ and 
$R_{23}(z)$ act in $C^n\otimes C^n\otimes C^n$ with
$R_{12}(z)=R(z)\otimes 1, R_{23}(z)=1\otimes R(z)$, ect. 
One can find that the R-matrix also satisfy
the following unitary and cross-unitary properties.
\begin{equation}
R_{12}(z_1-z_2)R_{21}(z_2-z_1)= \rho (z_1-z_2)\cdot id.
\end{equation}
\begin{equation}
R_{21}^{t_2}(z_2-z_1-nw)R_{12}^{t_2}(z_1-z_2)= 
\tilde{\rho}(z_1-z_2)\cdot id.
\end{equation}
where
\begin{equation}
\rho (z)=\frac{h(z+w)h(-z+w)}{h^2(w)}
\end{equation}
\begin{equation}
\tilde{\rho}(z)=\frac{h(z)h(-z-nw)}{h^2(w)},
\end{equation}
$w$ is defined by $w \equiv n\eta$ and 
$h(z)\equiv \sigma_0(z)$, $t_i$ means the transposition
in the $i$-th space. 

Assume an operator matrix L(z) satisfy the Yang-Baxter relation(YBR)
\begin{equation}
R_{12}(z_1-z_2)L_1(z_1)L_2(z_2)=L_2(z_2)L_1(z_1)R_{12}(z_1-z_2) ,
\end{equation}
with $L_1(z_1)=L(z_1)\otimes 1, L_2(z_2)=1\otimes L(Z_2)$. 
For the periodic boundary conditions,
its transfer matrix is defined by $t(z)=trL(z)$. 
We can prove that this tranfer matrices with different spectrum commute
with each other $[t(z_1),t(z_2)]=0$.
For the open boundary conditions, 
Sklyanin [23] proposed a systematic approach to  handle the problems which
involves the reflection equation (RE)
\begin{equation}
R_{12}(z_1-z_2)K_1(z_1)R_{21}(z_1+z_2)K_2(z_2)
=K_2(z_2) R_{12}(z_1+z_2)K_1(z_1)R_{21}(z_2-z_1), 
\end{equation}
where the reflecting $K$ matrix is a solution of the RE. 
In order to construct
the integrable models, we also need a dual reflection equation(DRE) 
which is associated with
the cross-unitary relation of the R matrix. 
For $Z_n$ belavin model, the dual RE takes the
form [24],
\begin{eqnarray}
& & R_{12}(z_1-z_2)\tilde{K}_1(z_1)R_{21}(-z_1-z_2-nw) \tilde{K}_2(z_2)
\nonumber\\
&=& \tilde{K}_2(z_2) R_{12}(-z_1-z_2-nw) \tilde{K}_1(z_1)R_{21}(z_2-z_1).
\end{eqnarray}
If we define the transfer matrix 
as $t(z)=tr[ \tilde{K}(z)L(z)K(z)L^{-1}(-z)] $,
with the help of unitary and cross-unitary relations, 
one can prove that $[t(z_1),t(z_2)]=0$, 
that means the model under consideration is integrable.

Additionally, there is an isomorphism between $K(z)$ and $\tilde{K}(z)$
\begin{equation}
\Phi: K(z)\rightarrow \tilde{K}(z)=K(-z-\frac{nw}{2}).
\label{11}
\end{equation}
Given a solution K(z) of the RE, we can find a solution 
$\tilde{K}(z)$ of the DRE.
For $Z_n$ Belavin model, a solution K(z) to the RE is as follows[24],
\begin{equation}
K(z)=K_0(z)K_0(0),
\label{12}
\end{equation}
where
\begin{equation}
K_0(z)=\sum_{\alpha \in Z^2_n}U_{2\alpha}(z)
\omega^{2\alpha_1\alpha_2}I_{2\alpha},
\end{equation}

\begin{equation}
U_{2\alpha}(z)=\frac{\sigma_{2\alpha}(z+c)}{\sigma_{2\alpha}(c)} ,
\end{equation}
$c$ is a arbitrary constant.

\section{The factorized L matrix}
Jimbo et al [25] defined the intertwiner of $Z_n$ model as an 
n-element column
vector $\phi_{a, a+\hat{\mu}}(z)$ whose $j$-th element is
\begin{equation}
\phi^{(j)}_{a, a+\hat{\mu}}(z)=\theta^{(j)}(z-nw\bar{a}_{\mu},n\tau),
\end{equation}
\begin{equation}
\theta^{(j)}(z-nw\bar{a}_{\mu},n\tau)=\theta \left [\begin{array}{c}
\frac{1}{2}-\frac{j}{n}\\
\frac{1}{2}
\end{array}\right ](z-nw\bar{a}_{\mu},n\tau),
\end{equation}
$\bar{a}_{\mu}=a_{\mu}-\frac{1}{n}\sum_{\nu}a_{\nu}+\delta_{\mu}, 
a\in Z^n$, $\delta_{\mu}$'s
are some generic numbers. Using intertiwner, 
the face-vertex correspondence can be
written as
\begin{eqnarray}
& & R_{12}(z_1-z_2)^{i'j'}_{ij}\phi^{(i')}_{a-\hat{\mu}-\hat{\nu}, 
a-\hat{\mu}}(z_1)\phi^{(j')}_{a-\hat{\mu}, a}(z_2) \nonumber\\
& = &\sum_{\kappa} W\left [\begin{array}{clcr}
a-\hat{\mu}-\hat{\nu} & a-\hat{\mu} \\
a-\hat{\kappa} & a
\end{array}\right ](z_1-z_2)\phi^{(j)}_{a-\hat{\mu}-\hat{\nu}, 
a-\hat{\kappa}}(z_2)
\phi^{(i)}_{a-\hat{\kappa}, a}(z_1),
\label{17}
\end{eqnarray}
where $W\left [\begin{array}{clcr}
 a-\hat{\mu}-\hat{\nu} & a-\hat{\mu} \\
 a-\hat{\kappa} & a
\end{array}\right ](z) $ is the face Boltzmann weight of 
$A^{(1)}_{n-1}$ IRF model [25]. It is defined as
 \begin{eqnarray}
 W\left [\begin{array}{clcr}
a-2\hat{\mu} & a-\hat{\mu} \\
a-\hat{\mu} & a
\end{array}\right ](z) &=& \frac {h(z+w)}{h(w)},\nonumber\\
 W\left [\begin{array}{clcr}
 a-\hat{\mu}-\hat{\nu} & a-\hat{\mu} \\
 a-\hat{\mu} & a
\end{array}\right ](z) &=& \frac{h(z-a_{\mu\nu}w)}{h(-a_{\mu\nu}w)},\\
 W\left [\begin{array}{clcr}
a-\hat{\mu}-\hat{\nu} & a-\hat{\nu} \\
a-\hat{\mu} & a
\end{array}\right ](z) &=& \frac{h(z)h((a_{\mu\nu}+1)w)}
{h(w)h(a_{\mu\nu}w)}.
\end{eqnarray}
The other face Boltzmann weights are defined as zeroes. 
Where $a_{\mu\nu}=\bar{a}_{\mu}-\bar{a}_{\nu}$.
At the same time, we can also find the 
n-element row vectors $\tilde{\phi}$ and $\bar{\phi}$ which satisfy
the following relations[20-22,25]
\begin{equation}
\tilde{\phi}^{(k)}_{a-\hat{\mu}, a}(z)\phi^{(k)}_{a-\hat{\nu}, a}(z)
=\delta_{\mu\nu},\end{equation}
\begin{equation}
\bar{\phi}^{(k)}_{a, a+\hat{\mu}}(z)
\phi^{(k)}_{a, a+\hat{\nu}}(z)=\delta_{\mu\nu}.\end{equation}
The above equation can also be written as
\begin{eqnarray}
\sum_{\mu}{\phi}_{a-\hat{\mu}, a}(z)\tilde{\phi}_{a-\hat{\mu}, a}(z)&=&I, \\
\label{22}
\sum_{\mu}{\phi}_{a, a+\hat{\mu}}(z)\bar{\phi}_{a, a+\hat{\mu}}(z)&=&I.
\label{23}
\end{eqnarray}\\
Using those results,the face-vertex correspondence can 
be written in other forms 
\begin{eqnarray}
& & \tilde{\phi}^{(i)}_{a-\hat{\mu}, a}(z_1)
R_{12}(z_1-z_2)^{i'j'}_{ij}\phi^{(j')}_{a-\hat{\nu}, a}(z_2) \nonumber\\
&=&\sum_{\kappa} W\left [\begin{array}{clcr}
a-\hat{\mu}-\hat{\kappa} & a-\hat{\nu} \\
a-\hat{\mu} & a
\end{array}\right ](z_1-z_2)\phi^{(j)}_{a-\hat{\mu}
-\hat{\kappa}, a-\hat{\mu}}(z_2)\tilde{\phi}^{(i')}_{a
-\hat{\mu}-\hat{\kappa}, a-\hat{\nu}}(z_1),
\end{eqnarray}

\begin{eqnarray}
& & \bar{\phi}^{(j)}_{a, a+\hat{\mu}}(z_2)
R_{12}(z_1-z_2)^{i'j'}_{ij}
\phi^{(i')}_{a-\hat{\nu}-\hat{\kappa}, a-\hat{\mu}}(z_1) \nonumber\\
& = &\sum_{\kappa} W\left [\begin{array}{clcr}
a & a+\hat{\nu} \\
a+\hat{\mu} & a+\hat{\mu}+\hat{\kappa}
\end{array}\right ](z_1-z_2)
\phi^{(i)}_{a+\hat{\mu}, a+\hat{\mu}+\hat{\kappa}}(z_1)
\bar{\phi}^{(j')}_{a+\hat{\nu}, a+\hat{\mu}+\hat{\kappa}}(z_2),
\end{eqnarray}

\begin{eqnarray}
& & \tilde{\phi}^{(i)}_{a-\hat{\mu}, a}(z_1)
\tilde{\phi}^{(j)}_{a-\hat{\mu}-\hat{\nu}, a-\hat{\mu}}(z_2)
R_{12}(z_1-z_2)^{i'j'}_{ij} \nonumber\\
& = &\sum_{\kappa} W\left [\begin{array}{clcr}
a-\hat{\mu}-\hat{\nu} & a-\hat{\kappa} \\
a-\hat{\mu} & a
\end{array}\right ](z_1-z_2)
\tilde{\phi}^{(j')}_{a-\hat{\kappa}, a}(z_2)
\tilde{\phi}^{(i')}_{a-\hat{\mu}-\hat{\nu}, a-\hat{\kappa}}(z_1),
\label{26}
\end{eqnarray}

\begin{eqnarray}
& & \bar{\phi}^{(i)}_{a+\hat{\mu}, a+
\hat{\mu}+\hat{\nu}}(z_1)\bar{\phi}^{(j)}_{a, a+\hat{\mu}}(z_2)
R_{12}(z_1-z_2)^{i'j'}_{ij} \nonumber\\
& = &\sum_{\kappa} W\left [\begin{array}{clcr}
a & a+\hat{\kappa} \\
a+\hat{\mu} & a+\hat{\mu}+\hat{\nu}
\end{array}\right ](z_1-z_2)
\bar{\phi}^{(j')}_{a+\hat{\kappa}, a+\hat{\mu}+\hat{\nu}}(z_2)
\bar{\phi}^{(i')}_{a, a+\hat{\kappa}}(z_1).
\end{eqnarray}
With the help of those face-vertex 
correspondence relations, we then can construct
the $L(z)$ matrix [20-22] which meet the YBR. Let
\begin{equation}
f(a, \mu, z)_i^j=\phi^{(i)}_{a-\hat{\mu}, a}(z+\xi_1)\tilde{\phi}^{(j)}_{a-
\hat{\mu}, a}(z+\xi_2)
\end{equation}
$\xi_1$ and $\xi_2$ are arbitrary complex numbers. 
Then from the face-vertex correspondence relations 
Eq.(\ref{17}) and (\ref{26}) , we have
\begin{eqnarray}
& & R(z_1-z_2)^{i'j'}_{ij} 
f(a-\hat{\mu}, \nu, z_1)_{i'}^{i''}f(a, \nu,z_2)_{j'}^{j''}+\mu 
\longleftrightarrow \nu\nonumber\\
&=&  f(a-\hat{\mu}, \nu,z_2)_{j}^{j'}
f(a, \nu, z_1)_{i}^{i'}R(z_1-z_2)^{i''j''}_{i'j'}+\mu 
\longleftrightarrow \nu ,
\end{eqnarray}
$\mu \longleftrightarrow \nu$ means the same form as the term before them, 
while the $\mu$ and $\nu$ exchange to each other.
We introduce the difference operator $\Gamma_{\mu}$
\begin{equation}
\Gamma_{\mu}f(a)=f(a+\hat{\mu})\Gamma_{\mu},\end{equation}
and we have
\begin{eqnarray}
& & R(z_1-z_2)^{i'j'}_{ij} 
\Gamma_{\nu}f(a, \nu, z_1)_{i'}^{i''}
\Gamma_{\mu}f(a, \nu, z_2)_{j'}^{j''}+\mu 
\longleftrightarrow \nu\nonumber\\
&=&  \Gamma_{\nu}f(a, \nu, z_2)_{j}^{j'}
\Gamma_{\mu}f(a, \nu, z_1)_{i}^{i'}R(z_1-z_2)^{i''j''}_{i'j'}+\mu 
\longleftrightarrow \nu .
\end{eqnarray}
We can find that the above equation is  
sufficient to prove the following YBR
\begin{equation}
R(z_1-z_2)^{i'j'}_{ij} L(a,z_1)_{i'}^{i''}L(a, z_2)_{j'}^{j''} =  
L(a, z_2)_{j}^{j'}L(a, z_1)_{i}^{i'}
R(z_1-z_2)^{i''j''}_{i'j'} .
\end{equation}
Where $L(a,z_k)=\sum_{\mu}\Gamma_{\mu}f(a,\mu, z_k)$. Let
\begin{equation}
g(a, \mu, z)_i^j=\phi^{(i)}_{a, a+\hat{\mu}}(z+\xi_2)
\bar{\phi}^{(j)}_{a, a+\hat{\mu}}(z+\xi_1)\end{equation}
There is
\begin{equation}
\sum_{\nu}\Gamma_{\nu} f(a, \nu, z)_i^j\sum_{\mu}
\Gamma_{-\mu}g(a, \mu, z)_j^k=\delta_{ik}\end{equation}
It shows that $L^{-1}(a, z)=\sum_{\mu}\Gamma_{-\mu}g(a, \mu, z)$ . 
Then the transfer matrix with open boundary conditions is
\begin{eqnarray}
t(z) &=& tr\tilde{K}(z)L(a, z)K(z)L^{-1}(a,-z)\nonumber\\
&=& tr\tilde{K}(z)\sum_{\nu}\Gamma_{\nu}f(a, \nu, z)K(z)\sum_{\mu}
\Gamma_{-\mu}g(a, \mu, -z)\nonumber\\
&=& \sum_{\mu \nu}\Gamma_{-\mu}\Gamma_{\nu}F^{(1)}_{\mu\nu}(a, z)
F^{(2)}_{\mu\nu}(a, z)\nonumber\\
&=& \sum_{\mu \nu}\Gamma_{-\mu}\Gamma_{\nu}G_{\mu\nu}(a, z),
\end{eqnarray}
with
\begin{equation}
F^{(1)}_{\mu\nu}(a,z)=\tilde{\phi}^{(i)}_{a+\hat{\mu}-\hat{\nu}, 
a+\hat{\mu}}(z+\xi_2)K(z)_i^j \phi^{(j)}_{a, a+\hat{\mu}}(-z+\xi_2)
\label{36}
\end{equation}
\begin{equation}
F^{(2)}_{\mu\nu}(a,z)=\bar{\phi}^{(k)}_{a, a+
\hat{\mu}}(-z+\xi_1)\tilde{K}(z)_k^l 
\phi^{(l)}_{a+\hat{\mu}-\hat{\nu}, a+\hat{\mu}}(z+\xi_1)
\label{37}
\end{equation}
\begin{equation}
G_{\mu\nu}(a,z)=F^{(1)}_{\mu\nu}(a, z)F^{(2)}_{\mu\nu}(a, z) .
\end{equation}

\section{The calculation of tansfer matrix t(z)}

Substituting the $K(z)$(\ref{12}) into Eq.(\ref{36}) , we get
\begin{equation}
F^{(1)}_{\mu\nu}(a,z)=\tilde{\phi}^{(i)}_{a+\hat{\mu}-
\hat{\nu}, a+\hat{\mu}}(z)\sum_{\gamma \in Z^2_n}
U_{2\gamma}(z)\omega^{2\gamma_1\gamma_2}[g^{2\gamma_2}
h^{2\gamma_1}K_0(0)]_i^j \phi^{(j)}_{a, a+\hat{\mu}}(-z).
\label{39}
\end{equation}
Here we let $\xi_1=\xi_2=0$. We can prove that
\begin{equation}
g^{\beta-1}h^{\alpha}K_0(0) \phi_{a, a+\hat{\mu}}(z)=(-1)^{\beta-1}
e^{2\pi\sqrt{-1}
\frac{\alpha}{n}(\frac{\alpha\tau}{2}+z-nw\bar{a}_{\mu}+\frac{1}{2})}
\phi_{a, a+\hat{\mu}}(-z+2nw\bar{a}_{\mu}-\alpha\tau-\beta),
\label{40}
\end{equation}
(see Apendix). 
The above Eq.(\ref{39}) can be written as
\begin{eqnarray}
F^{(1)}_{\mu\nu}(a,z) &=& \sum_{\gamma \in Z^2_n}U_{2\gamma}(z)
e^{2\pi\sqrt{-1}\frac{2\gamma_1\gamma_2}{n}}e^{2\pi\sqrt{-1}
\frac{2\gamma_1}{n}
(\frac{2\gamma_1\tau}{2}-z-nw\bar{a}_{\mu}+\frac{1}{2})}\nonumber\\
&\times& \tilde{\phi}^{(i)}_{a+\hat{\mu}-\hat{\nu}, a+\hat{\mu}}(z)
\phi_{a, a+\hat{\mu}}^{(i)}(z+2nw\bar{a}_{\mu}-2\gamma_1\tau-2\gamma_2-1)
\end{eqnarray}
From Eq.(\ref{22}) we know 
$\tilde{\phi}_{a+\hat{\mu}-\hat{\nu}, a+\hat{\mu}}(z)$ can be obtained
from the inverse of the matrix $\tilde{M}$ whose elements are
\begin{eqnarray}
\tilde{M}_{i\lambda}(z) &=& \phi^{(i)}_{a+\hat{\mu}-
\hat{\lambda}, a+\hat{\mu}}(z)\nonumber\\
&=& \theta^{(i)}(z-nw(\bar{a}_{\lambda}+\delta_{\mu\lambda}-1)).
\end{eqnarray}
So we get
\begin{equation}
\tilde{\phi}^{(i)}_{a+\hat{\mu}-\hat{\nu}, a+\hat{\mu}}(z)\phi_{a, a+
\hat{\mu}}^{(i)}(z+2nw\bar{a}_{\mu}-2\gamma_1\tau-(2\gamma_2+1))
=\frac{det\tilde{M}'}{det\tilde{M}}.
\end{equation}
Substituting the $\nu$ column elements of 
matrix $\tilde{M}$ with the corresponding elements
of column vector $\phi_{a, a+\hat{\mu}}(z
+2nw\bar{a}_{\mu}-2\gamma_1\tau-(2\gamma_2+1)) $ 
while keeping other matrix elements unchanged,
we get the matrix $\tilde{M}'$. 

Suppose the elments of a matrix A defined as 
$A_{ij}=\theta^{(i)}(nz_j)$, one
can prove that the determinant of the matrix A have the results [20-22]
 \begin{equation}
  detA=C(\tau)h(\sum_iz_i-\frac{n-1}{2})\prod_{i<k}h(z_i-z_k).
  \label{44}
  \end{equation}
  Using this result, we get
  \begin{eqnarray}
& & \tilde{\phi}^{(i)}_{a+\hat{\mu}-\hat{\nu}, a+\hat{\mu}}(z)
\phi_{a, a+\hat{\mu}}^{(i)}(z+2nw\bar{a}_{\mu}
-2\gamma_1\tau-(2\gamma_2+1)) \nonumber\\
&=& \frac{h(-z+w\delta+w(1-n)+\frac{n-1}{2}
-w(\bar{a}_{\nu}+\delta_{\mu\nu}-1+\bar{a}_{\mu})
+\frac{1}{n}(2\gamma_1\tau+2\gamma_2+1))}{h(-z
+w\delta+w(1-n)+\frac{n-1}{2})}\nonumber\\
&\times& \prod_{j\not=\nu}\frac{h(-w(\bar{a}_{j}
+\delta_{\mu j}-1+\bar{a}_{\mu})+\frac{1}{n}
(2\gamma_1\tau+2\gamma_2+1))}{h(-w(\bar{a}_{j}
+\delta_{\mu j}-\bar{a}_{\nu}-\delta_{\mu\nu}))},
\end{eqnarray}
($\delta=\sum_i \delta_i$). With the help of the formula
\begin{equation}
\sigma_0(z+\frac{1}{n}(\alpha\tau+\beta))=e^{-2\pi\sqrt{-1}
\frac{\alpha}{n}(\frac{\alpha\tau}{2n}+z+\frac{1}{2}
+\frac{\beta}{n})}\sigma_{\alpha,\beta}(z,\tau),
\label{46}
\end{equation}
finally, we get
\begin{eqnarray}
F^{(1)}_{\mu\nu}(a,z) &=& \sum_{\gamma \in Z^2_n}
U_{2\gamma}(z)e^{-2\pi\sqrt{-1}\frac{2\gamma_1\gamma_2}{n}}e^{2\pi\sqrt{-1}
\frac{2\gamma_1}{n}}\nonumber\\
&\times& \frac{\sigma_{2\gamma_1,2\gamma_2+1}(-z+w\delta+w(1-n)
+\frac{n-1}{2}-w(\bar{a}_{\nu}+\delta_{\mu\nu}-1+\bar{a}_{\mu}))}
{\sigma_0(-z+w\delta+w(1-n)+\frac{n-1}{2})}\nonumber\\
&\times& \prod_{j\not=\nu}\frac{\sigma_{2\gamma_1,2\gamma_2+1}
(-w(\bar{a}_{j}+\delta_{\mu j}-1+\bar{a}_{\mu}))}{\sigma_0
(-w(\bar{a}_{j}+\delta_{\mu j}-\bar{a}_{\nu}-\delta_{\mu\nu}))}.
\end{eqnarray}
Substituting the $\tilde{K}(z)$ into Eq.(\ref{37}) , considering the 
ismorphism relation between $\tilde{K}(z)$ and $K(z)$(\ref{11}), we get
\begin{equation}
F^{(2)}_{\mu\nu}(a,z)=\bar{\phi}^{(k)}_{a, a+\hat{\mu}}(-z)
\sum_{\gamma \in Z^2_n}U_{2\gamma}(-z-\frac{nw}{2})
\omega^{2\gamma_1\gamma_2}[g^{2\gamma_2}h^{2\gamma_1}
K_0(0)]_k^l \phi^{(l)}_{a+\hat{\mu}-\hat{\nu}, a
+\hat{\mu}}(z)\end{equation}
Using Eq(\ref{40}), we have
\begin{eqnarray}
F^{(2)}_{\mu\nu}(a,z) &=& \sum_{\gamma \in Z^2_n}U_{2\gamma}(-z
-\frac{nw}{2})e^{2\pi\sqrt{-1}\frac{2\gamma_1\gamma_2}{n}}
e^{2\pi\sqrt{-1}\frac{2\gamma_1}{n}
(\frac{2\gamma_1\tau}{2}+z-nw(\bar{a}_{\nu}
+\delta_{\mu\nu}-1)+\frac{1}{2})}\nonumber\\
&\times& \bar{\phi}^{(k)}_{a,a+\hat{\mu}}(-z)
\phi_{a+\hat{\mu}-\hat{\nu}, a+\hat{\mu}}^{(k)}
(-z+2nw(\bar{a}_{\nu}+\delta_{\mu\nu}-1)-2\gamma_1\tau-2\gamma_2-1)
\end{eqnarray}
From Eq.(\ref{23}), we also know $\bar{\phi}_{a, a+\hat{\mu}}(-z)$ 
can be obtained
from the inverse of the matrix $\bar{M}$ whose elements are
\begin{eqnarray}
\bar{M}_{i\lambda}(-z) &=& \phi^{(i)}_{a, a+\hat{\lambda}}(-z)\nonumber\\
&=& \theta^{(i)}(-z-nw\bar{a}_{\lambda}).
\end{eqnarray}
So we get
\begin{equation}
\bar{\phi}^{(k)}_{a, a+\hat{\mu}}(-z)
\phi_{a+\hat{\mu}-\hat{\nu}, a+\hat{\mu}}^{(k)}
(-z+2nw(\bar{a}_{\nu}+\delta_{\mu\nu}-1)-2\gamma_1\tau-2\gamma_2-1)
=\frac{det\bar{M}'}{det\bar{M}}.
\end{equation}
Substituting the $\mu$ column elements of matrix $\bar{M}$ 
with the corresponding elements
of column vector $\phi_{a+\hat{\mu}-\hat{\nu}, a+\hat{\mu}}
(-z+2nw(\bar{a}_{\nu}+\delta_{\mu\nu}-1)-2\gamma_1\tau-2\gamma_2-1) $, 
holding other matrix elements unchanged,
we get the matrix $\bar{M}'$. 
Using this result of (\ref{44}), we get

\begin{eqnarray}
& & \bar{\phi}^{(k)}_{a,a+\hat{\mu}}(-z)
\phi_{a+\hat{\mu}-\hat{\nu}, a+\hat{\mu}}^{(k)}(-z+2nw(\bar{a}_{\nu}
+\delta_{\mu\nu}-1)-2\gamma_1\tau-2\gamma_2-1)\nonumber\\
&=& \frac{h(z+w\delta+\frac{n-1}{2}-w(\bar{a}_{\mu}
+\delta_{\mu\nu}-1+\bar{a}_{\nu})+\frac{1}{n}
(2\gamma_1\tau+2\gamma_2+1))}{h(z+w\delta+\frac{n-1}{2})}\nonumber\\
&\times& \prod_{j\not=\mu}\frac{h(-w(\bar{a}_{j}
+\delta_{\mu \nu}-1+\bar{a}_{\nu})
+\frac{1}{n}(2\gamma_1\tau+2\gamma_2+1))}
{h(-w(\bar{a}_{j}-\bar{a}_{\mu}))}.
\end{eqnarray}

Similarly, with the help of Eq.(\ref{46}),
we get
\begin{eqnarray}
F^{(2)}_{\mu\nu}(a,z) &=& \sum_{\gamma \in Z^2_n}
U_{2\gamma}(-z-\frac{nw}{2})e^{-2\pi\sqrt{-1}
\frac{2\gamma_1\gamma_2}{n}}e^{2\pi\sqrt{-1}\frac{2\gamma_1}{n}}\nonumber\\
&\times& \frac{\sigma_{2\gamma_1,2\gamma_2+1}(z+w\delta
+\frac{n-1}{2}-w(\bar{a}_{\mu}+
\delta_{\mu\nu}-1+\bar{a}_{\nu}))}
{\sigma_0(z+w\delta+\frac{n-1}{2})}\nonumber\\
&\times& \prod_{j\not=\mu}\frac{\sigma_{2\gamma_1,
2\gamma_2+1}(-w(\bar{a}_{j}+\delta_{\mu\nu}-1
+\bar{a}_{\nu}))}{\sigma_0(-w(\bar{a}_{j}-\bar{a}_{\nu}))}.
\end{eqnarray}

\section{Taking trigonometric and scaling limit }

When $\tau \rightarrow \sqrt{-1}\infty$,  
considering $\frac{2\gamma_1}{n}\not=m$(m is integer number),
we can find that
\begin{equation}
\theta\left [\begin{array}{clcr}
\frac{1}{2}+\frac{2\gamma_1}{n}\\
\frac{1}{2}+\frac{2\gamma_2+1}{n}
\end{array}\right ](z,\tau) \rightarrow 
f(\tau)e^{2\pi\sqrt{-1}(\frac{2\gamma_1}{n}+
\frac{1}{2}+m)(z+\frac{1}{2}+\frac{2\gamma_2+1}{n})}.
\end{equation}
Then,
\begin{equation}
U_{2\gamma}(z)\rightarrow e^{2\pi\sqrt{-1}(\frac{2\gamma_1}
{n}+\frac{1}{2}+m)z}
\end{equation}
\begin{eqnarray}
F^{(1)}_{\mu\nu}(a,z) &\rightarrow& \sum_{\gamma_1}
\frac{e^{2\pi\sqrt{-1}(\frac{2\gamma_1}{n}+\frac{1}{2}+m)z}
e^{2\pi\sqrt{-1}\frac{2\gamma_1}{n}}}{\sin \pi(-z+w\delta
+w(1-n)+\frac{n-1}{2})} \nonumber\\
&\times& \prod_{j\not=\nu}\frac{1}{\sin 
\pi (\bar{a}_{j}+\delta_{\mu j}-\bar{a}_{\nu}
-\delta_{\mu\nu})} \nonumber\\
&\times& \sum_{\gamma_2}e^{-2\pi\sqrt{-1}
\frac{2\gamma_1\gamma_2}{n}}f^n(\tau)e^{2\pi
\sqrt{-1}(\frac{2\gamma_1}{n}+\frac{1}{2}+m)
(-z+n-\frac{1}{2}-nw\bar{a}_{\mu}+2\gamma_2)}\nonumber\\
&\leadsto& \sum_{\gamma_2}e^{2\pi\sqrt{-1}\frac{2\gamma_1
\gamma_2}{n}}=\sum_{\gamma_2}(e^{2\pi\sqrt{-1}
\frac{2\gamma_1}{n}})^{\gamma_2}=0.
\label{56}
\end{eqnarray}
We can see that from 
$F^{(2)}_{\mu\nu}(a,z)$, we can obtain the same result. 
So we only pay our attention to the case that
$\frac{2\gamma_1}{n} = m$. The next step is to consider this case.
  
We know there is an arbitrary parameter c 
in $U_{2\gamma}(z)$. Let $c=\epsilon\tau+c'$ with
  $\epsilon < \frac{1}{n}$, we have
   \begin{equation}
    \theta\left [\begin{array}{clcr}
      \frac{1}{2}+\frac{2\gamma_1}{n}\\
      \frac{1}{2}+\frac{2\gamma_2}{n}
  \end{array}\right ](z+c,\tau)=e^{-2\pi\sqrt{-1}
  \epsilon(\frac{2\epsilon\tau}{n}+z+c'+\frac{1}{2}+\frac{2\gamma_2}{n})}
    \theta\left [\begin{array}{clcr}
      \frac{1}{2}+\frac{2\gamma_1}{n}+\epsilon\\
      \frac{1}{2}+\frac{2\gamma_2}{n}
  \end{array}\right ](z+c',\tau),
  \end{equation}
   then
    \begin{eqnarray}
    U_{2\gamma}(z) &=& e^{-2\pi\sqrt{-1}\epsilon z}
    \frac{\theta\left [\begin{array}{clcr}
      \frac{1}{2}+\frac{2\gamma_1}{n}+\epsilon\\
      \frac{1}{2}+\frac{2\gamma_2}{n}
  \end{array}\right ](z+c',\tau)}{
    \theta\left [\begin{array}{clcr}
      \frac{1}{2}+\frac{2\gamma_1}{n}+\epsilon\\
      \frac{1}{2}+\frac{2\gamma_2}{n}
  \end{array}\right ](c',\tau)},\nonumber\\
  &\rightarrow& e^{-2\pi\sqrt{-1}
  \epsilon z}e^{2\pi\sqrt{-1}(\frac{1}{2}+\frac{2\gamma_1}{n}+\epsilon +m)z}
  \end{eqnarray}
  which have no relation to $\gamma_2$ and do not 
  include any dynamics variables like $\bar{a}_{\mu}$,
  so we handle it as an irrelevant constant when taking 
  trigonometic limit. The conclusion (\ref{56}) still 
  hold when $c=\epsilon\tau+c'$.
  Now, we can get$(\frac{2\gamma_1}{n}=m)$
    \begin{eqnarray}
  F^{(1)}_{\mu\nu}(a,z) &=& \sum_{\gamma_2}\frac{\sin 
  \pi(-z+w\delta+w(1-n)+\frac{n-1}{2}-w(\bar{a}_{\nu}+
  \delta_{\mu\nu}-1+\bar{a}_{\mu})+\frac{2\gamma_2+1}{n})}
  {\sin \pi(-z+w\delta+w(1-n)+\frac{n-1}{2})}\nonumber\\
 &\times& \prod_{j\not=\nu}\frac{\sin \pi(-w(\bar{a}_{j}+
 \delta_{\mu j}-1+\bar{a}_{\mu})+\frac{2\gamma_2+1}{n})}
 {\sin \pi(-w(\bar{a}_{j}+\delta_{\mu j}-\bar{a}_{\nu}-\delta_{\mu\nu}))},
 \label{59}
 \end{eqnarray}
   \begin{eqnarray}
  F^{(2)}_{\mu\nu}(a,z) &=& \sum_{\gamma_2}\frac{\sin \pi(z+w\delta+
  \frac{n-1}{2}-w(\bar{a}_{\nu}+\delta_{\mu\nu}-1+\bar{a}_{\mu})+
  \frac{2\gamma_2+1}{n})}{\sin \pi(z+w\delta+\frac{n-1}{2})}\nonumber\\
 &\times& \prod_{j\not=\mu}\frac{\sin \pi(-w(\bar{a}_{j}+\delta_{\mu\nu}
 -1+\bar{a}_{\mu})+\frac{2\gamma_2+1}{n})}{\sin \pi(-w(\bar{a}_{j}
 -\bar{a}_{\mu}))}.
 \label{60}
 \end{eqnarray}
 For the n conservations do not depend on the spectrum parameter z,
 when $z\rightarrow -\sqrt{-1}\infty$, we can get the integrable Hamiltonian 
 of the system.  So there are
   \begin{eqnarray}
  F^{(1)}_{\mu\nu}(a,z) &=& \sum_{\gamma_2}e^{\pi\sqrt{-1}(-w(\bar{a}_{\nu}
  +\delta_{\mu\nu}-1+\bar{a}_{\mu})+\frac{2\gamma_2+1}{n})}\nonumber\\
 &\times& \prod_{j\not=\nu}\frac{\sin \pi(-w(\bar{a}_{j}+\delta_{\mu j}-1
 +\bar{a}_{\mu})+\frac{2\gamma_2+1}{n})}{\sin \pi(-w(\bar{a}_{j}+\delta_{
 \mu j}-\bar{a}_{\nu}-\delta_{\mu\nu}))},
 \end{eqnarray}
   \begin{eqnarray}
  F^{(2)}_{\mu\nu}(a,z) &=& \sum_{\gamma_2}e^{-\pi\sqrt{-1}(-w(\bar{a}_{\nu}
  +\delta_{\mu\nu}-1+\bar{a}_{\mu})+\frac{2\gamma_2+1}{n})}\nonumber\\
 &\times& \prod_{j\not=\mu}\frac{\sin \pi(-w(\bar{a}_{j}
 +\delta_{\mu\nu}-1+\bar{a}_{\mu})+\frac{2\gamma_2+1}{n})}
 {\sin \pi(-w(\bar{a}_{j}-\bar{a}_{\mu}))}.
 \end{eqnarray}
 Expanding $\sin \pi(-w(\bar{a}_{j}+\delta_{\mu j}-1+\bar{a}_{\mu})
 +\frac{2\gamma_2+1}{n}) $ in $F^{(1)}_{\mu\nu}(a,z)$ and
 $\sin \pi(-w(\bar{a}_{j}+\delta_{\mu\nu}-1+\bar{a}_{\nu})
 +\frac{2\gamma_2'+1}{n}) $ in $F^{(2)}_{\mu\nu}(a,z)$ as the follows 
 resepectively,
 \begin{equation}
 \frac{1}{2\sqrt{-1}}[e^{\pi\sqrt{-1}(-w(\bar{a}_{j}
 +\delta_{\mu j}-1+\bar{a}_{\mu})+\frac{2\gamma_2+1}{n})} 
 -e^{-\pi\sqrt{-1}(-w(\bar{a}_{j}+\delta_{\mu j}-1+\bar{a}_{\mu})
 +\frac{2\gamma_2+1}{n})}],
 \end{equation}
 \begin{equation}
 \frac{1}{2\sqrt{-1}}[e^{\pi\sqrt{-1}(-w(\bar{a}_{j}+\delta_{\mu\nu}-1
 +\bar{a}_{\nu})+\frac{2\gamma_2'+1}{n})} - e^{-\pi\sqrt{-1}
 (-w(\bar{a}_{j}+\delta_{\mu\nu}-1+\bar{a}_{\nu})
 +\frac{2\gamma_2'+1}{n})}], 
 \end{equation}
 sum over $\gamma_2$ and $\gamma_2'$ and we can find that all 
 the terms of $F^{(1)}_{\mu\nu}(a,z) $ and $F^{(2)}_{\mu\nu}(a,z)$ have
 the forms
  \begin{equation}
  \sum_{\gamma_2}e^{\frac{2\pi\sqrt{-1}}{n}m\gamma_2} , m=1
  \underbrace{\pm 1\pm\cdots \pm 1}_{n-1} = n,n-2,n-4,\cdots,-n+2,
  \end{equation}
and
\begin{equation}
  \sum_{\gamma_2'}e^{\frac{2\pi\sqrt{-1}}{n}m\gamma_2'} , 
  m=1\underbrace{\pm 1\pm\cdots \pm 1}_{n-1} = n,n-2,n-4,\cdots,-n+2.
  \end{equation}
Only for case $ m=\pm n,0$,
the above summations are not equal to zeroes. 
When $n=even$,$m$ can take $n$ or 0, but it is rather complex to 
calculate the case $m=0$.
For simplicity, we only consider $n=odd$, 
so $m$ can only take $n$. We have
   \begin{equation}
  F^{(1)}_{\mu\nu}(a,z) \leadsto e^{\pi\sqrt{-1}
  \sum_i(-w(\bar{a}_{i}+\delta_{\mu i}-1+\bar{a}_{\mu})
  +\frac{1}{n})}\prod_{j\not=\nu}\frac{1}{\sin 
  pi(-w(\bar{a}_{j}+\delta_{\mu j}-\bar{a}_{\nu}-\delta_{\mu\nu}))},
  \label{67}
 \end{equation}
  \begin{equation}
  F^{(2)}_{\mu\nu}(a,z) \leadsto e^{-\pi\sqrt{-1}\sum_i(-w(\bar{a}_{i}
  +\delta_{\mu\nu}-1+\bar{a}_{\nu})+\frac{1}{n})}\prod_{j\not=\mu}
  \frac{1}{\sin \pi(-w(\bar{a}_{j}-\bar{a}_{\mu}))},
  \label{68}
 \end{equation}
  \begin{eqnarray}
  G_{\mu\nu}(a,z) &\leadsto& e^{\pi\sqrt{-1}(nw(\bar{a}_{\nu}
  -\bar{a}_{\mu})+nw\delta_{\mu\nu}-w)}\nonumber\\
  &\times& \prod_{j\not=\nu}\frac{1}{\sin \pi(-w(\bar{a}_{j}
  +\delta_{\mu j}-\bar{a}_{\nu}-\delta_{\mu\nu}))}\prod_{k\not=\mu}
  \frac{1}{\sin \pi(-w(\bar{a}_{k}-\bar{a}_{\mu}))}.
 \end{eqnarray}
For $\bar{a}_j-\bar{a}_k=a_j-a_k+\delta_j-\delta_k$,
and $H=\sum_{\mu\nu}\Gamma_{-\mu}\Gamma_{\nu}G_{\mu\nu}(a,z)$, 
in which $\Gamma_{\mu}f(a)=f(a+\hat{\mu})\Gamma_{\mu}$,
we get
  \begin{eqnarray}
  H &=& \sum_{\mu\nu}\Gamma_{\nu}[e^{\pi\sqrt{-1}(nw(a_{\nu}
  +\delta_{\nu}))}\prod_{j\not=\nu}\frac{1}{\sin \pi(-w(a_{j}
  -a_{\nu}+\delta_j-\delta_{\nu}))}]\nonumber\\
    &\times&  \Gamma_{-\mu}[e^{-\pi\sqrt{-1}(nw(a_{\mu}
    +\delta_{\mu}))}\prod_{k\not=\mu}\frac{1}{\sin \pi(-w(a_{k}
    -a_{\mu}+\delta_k-\delta_{\mu}))}]
 \end{eqnarray}
In quantum theory, one can set
$\hat{p}_{\mu}=\frac{\hbar}{\sqrt{-1}}\frac{\partial}{\partial x_{\mu}}$,
\begin{equation}
e^{\frac{\hbar}{\sqrt{-1}}\frac{\partial}
{\partial x_{\mu}}}f(x)=f(x+\frac{\hbar}{\sqrt{-1}}
\hat{\mu})e^{\frac{\hbar}{\sqrt{-1}}\frac{\partial}{\partial x_{\mu}}}.
\end{equation}
Compared with $\Gamma_{\mu}f(a)=f(a+\hat{\mu})\Gamma_{\mu}$, 
let $\Gamma_{\mu}\leadsto e^{\frac{\hbar}{\sqrt{-1}}
\frac{\partial}{\partial x_{\mu}}}$,
$-\sqrt{-1}w(a_{\mu}+\delta_{\mu}) \leadsto x_{\mu}$ and $w \leadsto \hbar$, 
there is

\begin{eqnarray}
H &\leadsto& \sum_{\mu\nu}e^{\hat{p}_{\nu}}[e^{-n\pi x_{\nu}}
\prod_{j\not=\nu}\frac{1}{\sinh \pi(x_{j}-x_{\nu})}]\nonumber\\
&\times& e^{-\hat{p}_{\mu}}[e^{n\pi x_{\mu}}\prod_{k\not=\mu}
\frac{1}{\sinh \pi(x_{k}-x_{\mu})}].
\label{72}
\end{eqnarray}
When $\hat{p}_{\lambda}\leadsto p_{\lambda}$, 
it becomes a Hamiltionan of a classical integrable syetem.

We also can begin directly from 
Eq.(\ref{59}) and (\ref{60}). Acorrding to the expanding method of getting
 the results (\ref{67}) and (\ref{68}), we have
   \begin{eqnarray}
  F^{(1)}_{\mu\nu}(a,z) &\leadsto& [e^{\pi\sqrt{-1}(-z+\frac{n-1}{2}+1-nw
  \bar{a}_{\mu})}-e^{-\pi\sqrt{-1}(-z+\frac{n-1}{2}+1-nw\bar{a}_{\mu}) }]
  \nonumber\\
 &\times& \prod_{j\not=\nu}\frac{1}{\sin \pi(-w(\bar{a}_{j}-\bar{a}_{\nu}
 +\delta_{\mu j}-\delta{\mu\nu}))},
 \end{eqnarray}
   \begin{eqnarray}
  F^{(2)}_{\mu\nu}(a,z) &\leadsto& [ e^{\pi\sqrt{-1}(z+\frac{n-1}{2}
  +1-nw(\bar{a}_{\nu}+\delta_{\mu\nu}-1))}-e^{-\pi\sqrt{-1}(z+
  \frac{n-1}{2}+1-nw(\bar{a}_{\nu}+\delta_{\mu\nu}-1))}]\nonumber\\
 &\times& \prod_{j\not=\mu}\frac{1}{\sin \pi(-w(\bar{a}_{j}
 -\bar{a}_{\mu}))},
 \end{eqnarray}
while ignoring those sin functions which have no any dynamics variables. Let
\begin{equation}
g_{\mu\nu}=\prod_{j\not=\nu}
\frac{1}{\sin \pi(-w(\bar{a}_{j}+\delta_{\mu j}-\bar{a}_{\nu}
-\delta_{\mu\nu}))}\prod_{k\not=\mu}\frac{1}{\sin \pi(-w(\bar{a}_{k}
-\bar{a}_{\mu}))},\end{equation}
we have
\begin{eqnarray}
 G_{\mu\nu}(a,z)&=&\{ e^{\pi\sqrt{-1}(n+1-nw(\bar{a}_{\mu}
 +\bar{a}_{\nu}+\delta_{\mu\nu}-1))}+e^{-\pi\sqrt{-1}(n+1-nw(\bar{a}_{\mu}
 +\bar{a}_{\nu}+\delta_{\mu\nu}-1))} \nonumber\\
 & & -e^{\pi\sqrt{-1}(-2z-nw(\bar{a}_{\mu}+1-\bar{a}_{\nu}
 -\delta_{\mu\nu}))}-e^{-\pi\sqrt{-1}(-2z-nw(\bar{a}_{\mu}+1
 -\bar{a}_{\nu}-\delta_{\mu\nu}))}\}g_{\mu\nu}.
 \end{eqnarray}
 For different $z$, $H(z)$ commute with each other, so we can obtain 
 three commutable
 Hamiltonians $H^i$ from $H$
 \begin{equation}
 H^i=\sum_{\mu\nu}\Gamma_{-\mu}\Gamma_{\nu}G^i_{\mu\nu}, i=1,2,3,
 \end{equation}
 with
 \begin{equation}
 G^1_{\mu\nu}= \{ e^{\pi\sqrt{-1}(n+1-nw(\bar{a}_{\mu}+\bar{a}_{\nu}
 +\delta_{\mu\nu}-1))}+e^{-\pi\sqrt{-1}(n+1-nw(\bar{a}_{\mu}
 +\bar{a}_{\nu}+\delta_{\mu\nu}-1))}\}g_{\mu\nu},\end{equation}
 \begin{equation}
 G^2_{\mu\nu}=  e^{\pi\sqrt{-1}(-2z-nw(\bar{a}_{\mu}+1-\bar{a}_{\nu}
 -\delta_{\mu\nu}))}g_{\mu\nu},
\end{equation}
\begin{equation}
G^3_{\mu\nu}=  e^{-\pi\sqrt{-1}(-2z-nw(\bar{a}_{\mu}+1-\bar{a}_{\nu}
-\delta_{\mu\nu}))}g_{\mu\nu}.
\end{equation}
In $G^1_{\mu\nu}$, if we have $\delta_{\mu}\rightarrow \delta_{\mu}
+\rho$, where $\rho $ is cummtable,
we may still get two integrable Hamiltonians
  \begin{eqnarray}
  H^{11} &=& \sum_{\mu\nu}\Gamma_{-\mu}\Gamma_{\nu}G^{11}_{\mu\nu}, \\
  H^{12} &=& \sum_{\mu\nu}\Gamma_{-\mu}\Gamma_{\nu}G^{12}_{\mu\nu},
  \end{eqnarray}
with
  \begin{eqnarray}
  G^{11}_{\mu\nu} &=& e^{\pi\sqrt{-1}(n+1-nw(\bar{a}_{\mu}
  +\bar{a}_{\nu}+\delta_{\mu\nu}-1))}g_{\mu\nu},\\
  G^{12}_{\mu\nu} &=& e^{-\pi\sqrt{-1}(n+1-nw(\bar{a}_{\mu}
  +\bar{a}_{\nu}+\delta_{\mu\nu}-1))}g_{\mu\nu}.
  \end{eqnarray}
$H^{11}$ and $H^{12}$ all can commute with $H^2$ and $H^3$, 
but they may not
commute with each other.

An integrable system require its poisson bracket
$H=f_1, f_2,\cdots, f_n$ satisfy
\begin{equation}
\{ f_i,f_j \}=\sum_{lm}(\frac{\partial f_i}{\partial p_l}
\frac{\partial f_j}{\partial q_m}-\frac{\partial f_i}{\partial q_l}
\frac{\partial f_j}{\partial p_m})=0.
\label{85}
 \end{equation}
When we take variables change, i.e. 
$q_i\rightarrow q_i'=\chi q_i+\rho_i $, $p_i$ remain unchange, 
in the new variable set $\{q_i',p_i\}$,
the above equation (\ref{85}) still hold, so the system is still integrable. 
Taking some kinds of limit after the variable subsitution, we may acheive
a new integrable system. i.e. in Eq.(\ref{72}), 
let $x_j\rightarrow x_j'=\chi x_j $, there is
\begin{eqnarray}
H' &=& \sum_{\mu\nu}e^{p_{\nu}}[e^{-n\pi \frac{1}{\chi}x_{\nu}'}
\prod_{j\not=\nu}\frac{1}{\sinh \frac{\pi}{\chi}(x_{j}'-x_{\nu}')}]
\nonumber\\
&\times& e^{-p_{\mu}}[e^{n\pi\frac{1}{\chi} x_{\mu}'}\prod_{k\not=\mu}
\frac{1}{\sinh \frac{\pi}{\chi}(x_{k}'-x_{\mu}')}].
\end{eqnarray}
When $\chi \rightarrow \infty$, take the lowest order term, we obtain
\begin{equation}
H'' = \sum_{\mu\nu}e^{p_{\nu}-p_{\mu}}\prod_{j\not=\nu}
\frac{1}{x_{j}-x_{\nu}}\prod_{k\not=\mu}\frac{1}{x_{k}-x_{\mu}}.
\end{equation}
This may be an integrable Hamiltonian.
 
Additionally, we can let $ x_0\rightarrow x_0+\rho_0$ 
with others unchanged in H, and let
$\sinh(\pi(x_0-x_i))$ become $\frac{1}{2}e^{\pi(x_0-x_i)}$ 
to obtain a new Hamiltonian. We also can change two variables or more 
to get a new Hamiltonian.
 i.e. let $\delta_0 \gg \delta_1 \gg \cdots \gg \delta_{n-1} 
 \gg 1\rightarrow \infty$,
 $x_i\rightarrow x_i+\delta_i$, there are $\sinh\pi(x_i-x_j) 
\rightarrow\frac{1}{2}e^{\pi(x_i-x_j +\delta_i-\delta_j))}$($i<j$). 
Then
\begin{eqnarray}
H &\rightarrow& \sum_{\mu\nu}e^{p_{\nu}-p_{\mu}}e^{-n\pi (x_{\mu}
-x_{\nu})}\nonumber\\
&\times& \prod_{j<\nu}e^{-\pi(x_j-x_{\nu})}\prod_{j>\nu}
e^{-\pi(x_{\nu}-x_j)}\prod_{k<\mu}e^{-\pi(x_k-x_{\nu})}
\prod_{k>\mu}e^{-\pi(x_{\nu}-x_k)}.
\end{eqnarray}
Like this, we can get amount of integrable Hamiltonians.

\section{Summary and discussions}
We study the $Z_n$ Belavin model with the open
boundary conditions. By the factorized $L$ operators, we
constructed the commuting difference operators which
is similar to the Ruijsenaars-Macdonald's difference operators.
We obtained a one-dimensional many-body integrable model from
this commuting difference operators.
By taking the special limit, this model is similar as the 
Calogero model.

The one-dimensional integrable models presented in this paper is obtained
from the transfer matrix operators with open boundary conditions.
Generally, for different solutions of the reflection equations, 
we can define different transfer matrices. Thus, to some extent, 
we give a systematic approach to obtain one-dimensional many-body
integrable models. Besides the $Z_n$ Belavin model, 
there are a lot of two-dimensional exactly
solvable model in statistical mechanics, it is interesting to
find relations between those models with some one-dimensional
many-body systems.

\vskip 2truecm
{\large Acknowlegements:} This work is supported in part
by the Natural Science Foundation of China.
\newpage

\appendix{\bf \Large Appendix }

\begin{eqnarray}
& & [g^{\beta}h^{\alpha}K_0(0) \phi_{a, a+\hat{\mu}}(z)]^{(i)}\nonumber\\
&=& \omega^{\beta i}\delta_{ij}\delta_{j+\alpha, k}
\bar{\delta}_{k+l, 0}\theta^{(l)}(z-nw\bar{a}_{\mu},n\tau)\nonumber\\
&=& \omega^{\beta i}\theta\left[\begin{array}{cl}
\frac{1}{2}+\frac{i+\alpha}{n}\\ \frac{1}{2}\end{array}
\right](z-nw\bar{a}_{\mu},n\tau)\nonumber\\
&=& \omega^{\beta i}\theta\left[\begin{array}{cl}-
\frac{1}{2}-\frac{i+\alpha}{n}\\ -\frac{1}{2}\end{array}
\right](-z+nw\bar{a}_{\mu},n\tau)\nonumber\\
&=& (-1)\omega^{\beta i}e^{2\pi\sqrt{-1}
\frac{i+\alpha}{n}}\theta\left[\begin{array}{cl}
\frac{1}{2}-\frac{i+\alpha}{n}\\ \frac{1}{2}
\end{array}\right](-z+nw\bar{a}_{\mu},n\tau)\nonumber.
\end{eqnarray}
Using the formula
\begin{equation}
\theta^{(i)}(z+\alpha\tau+\beta,n\tau)=e^{-2\pi\sqrt{-1}
\frac{\alpha}{2}(\frac{\alpha\tau}{2}+z+\frac{1}{2}+\beta)}
e^{2\pi(\frac{1}{2}-\frac{i}{n})\beta}\theta\left[\begin{array}
{cl}\frac{1}{2}-\frac{i-\alpha}{n}\\ \frac{1}{2}\end{array}
\right](z,n\tau)\end{equation}
there is
 \begin{eqnarray}
\cdots &=& (-1)^{(\beta-1)}e^{2\pi\sqrt{-1}\frac{\alpha}{n}
(\frac{\alpha\tau}{2}+z-nw\bar{a}_{\mu}+\frac{1}{2})}
\omega^{i}\theta\left[\begin{array}{cl}\frac{1}{2}-\frac{i}{n}
\\ \frac{1}{2}\end{array}\right]
(-z+nw\bar{a}_{\mu}-\alpha\tau-\beta,n\tau)\nonumber\\
&=& (-1)^{(\beta-1)}e^{2\pi\sqrt{-1}\frac{\alpha}{n}
\frac{\alpha\tau}{2}+z-nw\bar{a}_{\mu}+\frac{1}{2})}
\omega^{i}\phi^{(i)}_{a,a+\hat{\mu}}(-z+2nw\bar{a}_{\mu}
-\alpha\tau-\beta,n\tau)\nonumber\\
&=& (-1)^{(\beta-1)}e^{2\pi\sqrt{-1}\frac{\alpha}{n}
(\frac{\alpha\tau}{2}+z-nw\bar{a}_{\mu}+\frac{1}{2})}
[g\phi_{a,a+\hat{\mu}}(-z+2nw\bar{a}_{\mu}-\alpha\tau-\beta,n\tau)]^{(i)}
\end{eqnarray}
move the $g$ in the right hand of the equation to the left side 
and we get the conclusion (\ref{40}).

\end{document}